\documentclass[12pt]{article}

\usepackage[T2A]{fontenc}
\usepackage[utf8]{inputenc}

\usepackage{amsfonts}
\usepackage{amssymb}
\usepackage{amsmath}
\usepackage{amsthm}

\def \ge {\geqslant}

\theoremstyle{plain}

\usepackage{graphicx}
\begin{document}

\begin{Large}
\centerline{\bf
Simultaneous two-dimensional best Diophantine  }
\centerline{\bf approximations in the Euclidean norm}

\vskip+0.5cm
\centerline{Evgeny V. Ermakov}\begin{flushleft}
\begin{flushleft}

\end{flushleft}
\end{flushleft}
\vskip+1.0cm
\end{Large}

{\bf\Large 1. Introduction}\\

This paper is devoted to the exponents of growth of denominators of best simultaneous Diophantine approximations. Consider $\mathbb{R}^n$ with a norm $\|\cdot\|$. For any  vector $\alpha = (\alpha_1,\dots,\alpha_n) \in\mathbb{R}^n \backslash\mathbb{Q}^n $ and any $q\in\mathbb{Z}$ define the following value: $$\delta_q = \min\limits_{\textbf{p}=(p_1,\dots,p_n)\in\mathbb{Z}^n} \|q\cdot\alpha - \textbf{p}\|.$$
Let $\textbf{p}(q)\in \mathbb{Z}^n$  be the vector, where the minimum is attained; let $\textbf{r}(q)=q\cdot\alpha - \textbf{p}(q)$, so $\delta_q =\|\textbf{r}(q)\|.$ Given a norm $\|\cdot\|$ in $\mathbb{R}^n$ and a vector $\alpha\in\mathbb{R}^n \backslash\mathbb{Q}^n$ we can define the sequence of best approximations (with respect to this norm) as a sequence $(q_k)_{k=1}^\infty$, such that $q_1 =1$ and $\forall q<q_k \,\,\,\,\,\,\, \delta_q>\delta_{q_k}.$
Now we can define following values:
$$g(\alpha,\|\cdot\|)=\liminf\limits_{k\to\infty} (q_k)^{\frac{1}{k}} ,$$
$$G(n,\|\cdot\|)=\inf\limits_{\alpha\in\mathbb{R}^n \backslash\mathbb{Q}^n}g(\alpha ,\|\cdot\|) .$$
J.Lagarias \cite{1} has proved the following statement:
\\

{\bf Theorem 1.}\,\,\,
\emph{For any norm $\|\cdot\|$ on $\mathbb{R}^n$ and a vector $\alpha$, that has at least one irrational coordinate, the inequality $q_{k+2^{n+1}}\geq 2q_{k+1}+q_k$ holds for all  $k\geq 1.$ So $G(n,\|\cdot\|)\geq\theta$, where $\theta$ is the maximal positive root of $\theta^{2^{n+1}}=2\theta +1$.}
\\

In this paper we consider $\mathbb{R}^2$ with the Euclidian norm. From Theorem 1 it follows that for the Euclidian norm in $\mathbb{R}^2$, and any vector $\alpha$  one has
$q_{k+8} \geq q_{k+1} + q_k .$
\\

There is another well known statement that holds for any norm. Given a norm $\|\cdot\|$ in $\mathbb{R}^n$ consider the contact number  $K(n,\|\cdot\|)$. This number is defined as the maximal number of unit balls with respect to the norm$\|\cdot\|$ without interior common points that can touch another unit ball.
\\

{\bf Theorem 2.}\,\,\,
\emph{For any norm $\|\cdot\|$ on $\mathbb{R}^n$ with the contact number $K=K(n,\|\cdot\|)$ and a vector $\alpha$, that has at least one irrational coordinate, we have the inequality $q_{k+K}\geq q_{k+1}+q_k$, and so $G(\|\cdot\|)\geq\theta$, where $\theta$ is maximum positive root of $\theta^K=\theta +1$.}
\\

For the Euclidian norm in $\mathbb{R}^2$ we have $K=6$. So Theorem 2 gives  the inequality \begin{equation}\label{1.1}q_{k+6} \geq q_{k+1} + q_k.\end{equation} It follows that $G(2,\|\cdot\|_e)\geq\theta,$ where $\theta$ is maximum positive root of $\theta^6=\theta +1$ and $\|\cdot\|_e$ is the Euclidian norm.\\

Theorem 2 is a well known result, one can find a proof of it in M.Romanov paper \cite{3}. M.Romanov \cite{3} proved a stronger result that the inequality \begin{equation}\label{1.2}q_{k+4} \geq q_{k+1} + q_k.\end{equation} is valid for any $k\geq 1.$ From inequality (\ref{1.2}) it follows that $G(2,\|\cdot\|_e)\geq\theta_0$ where $\theta_0$ is a positive root of $\theta_0^6=\theta_0 +1, \,\,\,\, \theta_0 =1.220744$ The main result of the present paper is an improvement of Romanov's result.\\

{\bf Theorem 3.}\,\,\, \emph{For  the Euclidian norm
in $\mathbb{R}^2$ and any vector $\alpha$, that has at least one
irrational coordinate one has $G(2,\|\cdot\|_e)\geq 1.228043.$}
\\

The proof of Theorem 3 is based on following geometric statement that together with the inequality (\ref{1.2}) and some numerical calculations gives the lower bound.

{\bf Theorem 4.}\,\,\,
\emph{Suppose that  $\alpha\in\mathbb{R}^2$ has at least one irrational coordinate. Let $q_k \ldots q_{k+4}$ be consecutive denominators from the sequence of best approximations in Euclidian norm for vector $\alpha$. Then for every $k\geq 1$ at least one of two following inequalities are valid:}
\begin{equation}\label{2.1}q_{k+3} + q_{k+2}\geq 2q_{k+1} + q_k\end{equation}
\begin{equation}\label{2.2}q_{k+4}\geq q_{k+2} + q_k\end{equation}
\emph{Moreover, among any two successive values of $k$ for at least one value the inequality (\ref{2.1}) holds.}
\\

A.Brentjes \cite{2} gave the following example. Let $\eta$ be the maximal root of the equation $\eta^3=\eta+1, \,\,\,\,\,\,\,\eta= 1.3248...$. Then for $\alpha=(\alpha_1,\alpha_2)=(\eta,\eta^2)$
one has $g(\alpha,\|\cdot\|_e)=\eta$. J.Lagarias \cite{1} made a conjecture, that $G(2,(\|\cdot\|_e))=\eta$.

In Sections
 2,3 below we give a complete proof of Theorem 4. In Section 4 we deduce Theorem 3 from Theorem 4 and Romanov's theorem. There we describe all necessary computer calculations.\\
\\

{\bf\Large 2. Geometric lemmas} \\

{\bf Lemma 1.}\,\,\, \emph{Consider a convex hexagon $A_1 A_2 A_3 A_4 A_5 A_6$. Suppose that its opposite sides
are equal and parallel.
Suppose that
$O$ is an interior point of the  hexagon. Let
all the distances
  $|A_1O |, |A_2O |, |A_3O |, |A_5O |$ are different.
Then there exists $i \in \{ 1,2,3,5\}$ such that
$$|A_iO|> \min\limits_{j=1,2,3} {|A_j A_{j+1} |}.$$
} \\

\emph{Proof.}
Let $a= \min(|A_1 A_2|,|A_2 A_3|)$. Without loss of generality  suppose that $a= |A_1 A_2|$.

Consider circles $\omega_1$ and $\omega_2$ with radiuses $a$ and centers in $A_1$ and $A_3$ correspondingly. Let $\kappa_1$ and  $\kappa_2$ be closed disks bounded by $\omega_1$ and $\omega_2$. Define $\Omega = \kappa_1 \cap\kappa_2$. (See fig.1.)

Suppose that the conclusion of Lemma 1 is not true, that is  there
exists an  interior point O of hexagon $A_1 A_2 A_3 A_4 A_5 A_6$
such that $|A_1 O |, |A_2 O |, |A_3 O |, |A_5 O |$ are different
and
$$|A_iO|\leq \min\limits_{j=1,2,3} {|A_j A_{j+1} |}, \,\,\, i=1,2,3,5.$$

\begin{center}
\includegraphics[bb= 0 0 12cm 13cm]{s.1}
\end{center}
\setlength{\unitlength}{1cm}
\begin{picture}(50,.50)
\put(13,.5){fig.1}
\end{picture}

So there exist $i\in \{1,2,3,5\}$ such that
$$\|OA_i\|\leq \min\limits_{j=1,2,3} {|A_j A_{j+1} |}\leq a.$$
 By the  condition $\max(|OA_1|,|OA_3|)\leq a$ we see that $O\in\Omega$.
So $\Omega\neq\varnothing$ and circles $\omega_1$ and $\omega_2$ have common points. If $\omega_1$ and $\omega_2$ have
the  unique common point $O$ then $\|A_1O\| =\|A_3O\|$.  This contradicts to the conditions of Lemma 1. So we see that circles $\omega_1$ and $\omega_2$ have two different common points.

The line $A_1A_3$ divides the plane into two different half-planes. Define $Q$ to be that point of the  intersection $\omega_1$ and $\omega_2$ such that $A_2$ and $Q$ belong to different half-planes. Let $M$ be the point  symmetric to  $A_2$ with respect to the  center of the segment  $A_1A_3$. So $MA_3 A_4 A_5$ is a parallelogram and $M\in\omega_2.$ Consider the disk $\Theta$ with center in $A_5$ and radius $\|A_5M\|=\|A_3A_4\|$.\\

By the construction $O\in \Omega\cap\Theta.$ But if $\Omega$ and $\Theta$ have  a common point, it is the unique point $Q=M$ as the distance from $Q$ to the line $A_1A_3$ is less or equal to the distance from $M$ to the line $A_1A_3$. So $M$ belongs to $\omega_2$ but does not belong to $\Omega$ if it is not point of intersection of $\omega_1$ and $\omega_2.$ So if such point $O$ exists it is equal to $Q$. This  contradicts to the condition that $\|A_1O\|\neq \|A_3O\|.$ Lemma 1 is proved.\\

Suppose that $q_{k+3}< q_{k+1} + q_k$, otherwise we at once get (\ref{2.1}) as the sequence $(q_k)$ increases.

Consider remainder vectors $\textbf{r}(q_k), \textbf{r}(q_{k+1}), \textbf{r}(q_{k+2}), \textbf{r}(q_{k+3}).$ There exist a substitution of four indices $s=(s(1),s(2),s(3),s(4))$ such  that $\textbf{r}(q_{k-1+i})=\overrightarrow{OR_{s(i)}}$ and $R_1R_2R_3R_4$ is a tetragon without self intersections.\\

{\bf Lemma 2.}\,\,\, \emph{1. The tetragon $R_1R_2R_3R_4$ is convex, point O lies inside it.}

\emph{2. All of its sides and diagonals are not less then the longest remainder vector $|r(q_k)|$.}

\emph{3. Angles between vectors $\overline{OR_i}$ and $\overline{OR_j}$ ($i\neq j$) are greater than} $\frac{\pi}{3}$.\\

\emph{Proof.}
Suppose, that $|R_iR_j|<|r(q_k)|$ for any $i\neq j$. Let $R_i,R_j$ are the endpoints of
vectors $\textbf{r}(q_s)$ and $\textbf{r}(q_l)$ correspondingly. Then $|\textbf{r}(|q_s-q_l|)|<|\textbf{r}(q_k)|.$
From $q_{k+3}< q_{k+1} + q_k$ it follows that $0<|q_s-q_l|<q_{k+1}.$ Last inequalities contradict to the fact
that $q_k$ and $q_{k+1}$ are denominators of consecutive best approximations.
The second statement of Lemma 2 is proved.

In any triangle $OR_iR_j , i \neq j $ the side $ R_iR_j$ is the greatest one. Lengths of $\textbf{r}(q_k)$
decrease strictly, so those triangles can not have three equal sides and angles between
vectors $\overrightarrow{OR_i}$ are greater then $\frac{\pi}{3}.$
Other angles in these triangles are less or equal  to $\frac{\pi}{3}$.
We see that  $R_1R_2R_3R_4$ is convex, and the point O lies inside it. Lemma 2 is proved.\\
\\

{\bf\Large 3. Proof of Theorem 4}\\

We need two more lemmas.\\

{\bf Lemma 3.}\,\,\, \emph{If tetragon $R_1R_2R_3R_4$  is not a parallelogram, then the  inequality (\ref{2.1}) holds.}
\\

\emph{Proof.}
If $R_1R_2R_3R_4$ has no parallel sides, then we can make a convex hexagon by building parallelograms on two pairs of its sides.
(See fig.2.)
 Without loss of generality we may suppose that  the  hexagon vertex $R_4$ lies between the vertices $X_1$ and $X_2$. So we have constructed the hexagon $R_1R_2R_3X_2R_4X_1$.

\begin{center}
\includegraphics{s.2}
\end{center}
\setlength{\unitlength}{1cm}
\begin{picture}(50,.50)
\put(13,.5){fig.2}
\end{picture}

Consider the segment $R_3X_1$ (it is equal and parallel to segment $R_1X_2$).  Put $x=|R_3X_1|$.
By  the construction the length of the remainder vector for the denominator $q=|q_1+q_3-q_2-q_4|$ is not greater then $x$.
\\
As the sequence $(q_k)$ increases strictly, we have three possible values of $q$. So we should  consider three cases.

\textbf{Case 1.} \,\,\,  $q=|q_{k+3}+q_k-q_{k+2}-q_{k+1}|$. \,\,\, Here $0<q<q_k$, and the length of  the remainder vector for $q$ is not less then $|\textbf{r}(q_{k-1})|$. So $x\geq |\textbf{r}(q_k)|$.

\textbf{Case 2.} \,\,\,  $q=q_{k+3}+q_{k+1}-q_{k+2}-q_k$. \,\,\, Here $0<q<q_{k+1}$. The  length of the  remainder vector for $q$ is not less then $|\textbf{r}(q_k)|$ ($q$  is the  denominator of the next best approximation). So $x\geq |\textbf{r}(q_k)|$.

\textbf{Case 3.} \,\,\,  $q=q_{k+3}+q_{k+2}-q_{k+1}-q_k$. \,\,\, Then $q>0$ and we have 2 subcases:

\textbf{3a}  $q=q_{k+3}+q_{k+2}-q_{k+1}-q_k<q_{k+1}$. \,\,\, Here as in cases 1 and 2 we have  $x\geq |\textbf{r}(q_k)|$.

\textbf{3b}  $q=q_{k+3}+q_{k+2}-q_{k+1}-q_k\geq q_{k+1}$ \,\, Here we get the inequality (\ref{2.1}).\\

In  cases {\bf 1}, {\bf  2}, {\bf 3a} we have the following situation.
As $|\textbf{r}(q_k)|>|\textbf{r}(q_{k+1})|>|\textbf{r}(q_{k+2})|>|\textbf{r}(q_{k+3})|$ we see that  the hexagon
$R_1R_2R_3X_2R_4X_1$ and the zero point $O$
satisfy the conditions of Lemma 1.
By Lemma 1 we see that
$$\max_{i=1,2,3,4} |R_iO|> \min \{|X_1R_1|,|R_1R_2|,|R_2R_3|\}.$$
As in our  cases $x\geq |\textbf{r}(q_k)|$ we see that
$$\max_{i =1,2,3,4} |R_iO| > \min \{|R_4R_1|,|R_1R_2|,|R_2R_3|,|R_3R_4|\}.$$

This contradicts to Lemma 2.
So the cases are  {\bf 1}, {\bf  2}, {\bf 3a}  are not possible.

But in the remaining case
 {\bf 3b}  we have the  inequality (\ref{2.1}).
\\

To finish the proof of Lemma 3 we must consider the case when
$R_1 R_2 R_3 R_4$ has a pair of parallel sides.
Then the hexagon $R_1R_2R_3X_2R_4X_1$
is a degenerate one (two its angles are equal to $\pi$).
Now the proof follows the steps of the proof in non-degenerate case.
The only difference is that we
apply Lemma 1
for the  degenerate hexagon. Lemma 3 is proved.
\\

{\bf Lemma 4.}\,\,\, \emph{If  $R_1R_2R_3R_4$  is  a parallelogram and $q_{k+3}< q_{k+1} + q_k$, then endpoints of
the next four remainder vectors (for $k+1, k+2, k+3, k+4$) do not form a parallelogram.}

\emph{Proof.}
Suppose they do. Let $\textbf{r}(q_{k+4})= \overline{OR_5}, \,\,\textbf{r}(q_k)= \overline{OR}.$ This parallelogram has three common vertices with $R_1R_2R_3R_4$. So one of the vertices of the hexagon  $R_1R_2R_3R_4$ is the  center of the segment $RR_5.$
This vertex we denote by $R_6$.

As
$\|OR\|=|\textbf{r}(q_k)|>|\textbf{r}(q_{k+4})|=\|OR_5\|$,
we see that the zero point $O$ lies closer to  $R_5$ than to  $R.$ So in  the triangle $ORR_6$ the angle in the vertex  $R_6$ is greater than $\frac{\pi}{3}$ and  the  length of the remainder vector $\textbf{r}(q_k)= \overline{OR}$ is
greater than the length of the  parallelogram`s side $RR_2$.
We get the contradiction to
Lemma 2. Lemma 4 is proved.\\

\emph{Proof of Theorem 4.}

{\bf 1.}
If points $R_1, R_2, R_3, R_4$ do not form a parallelogram, then using Lemma 3 we get inequality (\ref{2.1}).

{\bf 2.} If the inequality $q_{k+3}< q_{k+1} + q_k$ do not holds, we again  get inequality (\ref{2.1}).

{\bf 3.} We may suppose that
$R_1, R_2, R_3, R_4$ do  form a parallelogram and $q_{k+3}< q_{k+1} + q_k$.
Then by Lemma 4 the  endpoints of  the next four remainder vectors (for $k+1, k+2, k+3, k+4$) do not form a parallelogram.
So for the  approximations $k+1, k+2, k+3, k+4$  the  inequality  (\ref{2.1})  is valid. We see that
\begin{equation}\label{2.3}q_{k+4} + q_{k+3}\geq 2q_{k+2} + q_{k+1}.\end{equation}

Let the endpoints of vectors
$$\textbf{r}(q_k), \textbf{r}(q_{k+1}), \textbf{r}(q_{k+2}), \textbf{r}(q_{k+3})$$
 form a parallelogram in
the   order
$$\textbf{r}(\hat{q}_1), \textbf{r}(\hat{q}_2), \textbf{r}(\hat{q}_3), \textbf{r}(\hat{q}_4).$$
 Then the  remainder vector for the denominator $p=|\hat{q}_1+\hat{q}_3-\hat{q}_2-\hat{q}_4|$
is equal to zero. As   $\alpha$ is not a rational vector we see that  $p=0$. As the sequence of denominators of best approximations increases strictly we get $0=p=q_k+q_{k+3}-q_{k+2}-q_{k+1}$.

The last equality  together with (\ref{2.3}) implies  (\ref{2.2}). Theorem 4 is proved.\\
\\

{\bf\Large 4. Proof of Theorem 3}\\

From Theorem 4  we immediately obtain

{\bf Proposition 1.}\,\,\, Let $ l \in \mathbb{R}$.
\emph{Let $\alpha\in\mathbb{R}^2 \backslash\mathbb{Q}^2.$ Let $\in\mathbb{R}$. Then for every $k\geq 1$ for
five consecutive denominators
$q_k ,\ldots , q_{k+4}$  we have at least one of  three following inequalities}
\begin{equation}\label{3.1}q_{k+2}\geq lq_{k+1}\end{equation}
\begin{equation}\label{3.2}q_{k+3}\geq (2-l)q_{k+1} + q_k\end{equation}
\begin{equation}\label{2.2}q_{k+4}\geq q_{k+2} + q_k\end{equation}
\emph{Moreover, for any two successive values of $k$ for at least one value  the
inequality (\ref{3.1}) or the inequality (\ref{3.2}) holds.}\\
\\

For further proof we need to use some computer calculations.

Let $0<l<2$. Let $ m=1, \ldots ,7 $.

Put $r_0=r_1=31,\,\,\,\, r_2=r_3=r_4=33,\,\,\,\, r_5=34,\,\,\,\,
r_6=35,\,\,\,\,$

$l_0=\ldots =l_3=1.298,\,\,\,\, l_4=l_5=l_6=1.293,\,\,\,\,$

$\theta_0=1.2207, \,\,\,\theta_1=1.2272, \,\,\,\theta_2=1.2275,
\,\,\,\theta_3=1.22779,$

$\theta_4=1.2278,\,\,\,\, \,\,\,\theta_5=1.22785,
\,\,\,\theta_6=1.22791, \,\,\,\,\theta_7=1.228043.$
\\

Consider a sequence $I=(i_0, \ldots ,i_{r-1}), \,\,\,\,
i_\nu\in\{1,2,3\}$ such that in any couple $i_\nu, i_{\nu+1}$ at
least one element is not equal to 3.
For such $I$ we construct a sequence $\{Q_k(I, m)\}, \,\,\,\, 0\leq k\leq r+3$ by
the following procedure.

First of all we define three rules for obtaining the vector
$$(Q^{j+1}_{j+1},\,
Q^{j+1}_{j+2},\, Q^{j+1}_{j+3},\, Q^{j+1}_{j+4})$$
 from the vector
$$(Q^{j}_{j},\, Q^{j}_{j+1},\, Q^{j}_{j+2},\, Q^{j}_{j+3}) :$$ rule
$\Re_1$, rule $\Re_2$ and rule $\Re_3$. These rules correspond to different inequalities in Proposition 1.\\

Rule $\Re_1$:

$$
\begin{cases}
 Q^{j+1}_{j+1} = Q^{j}_{j+1}, \cr
Q^{j+1}_{j+2} = \max\{lQ^{j}_{j+1},\,\, Q^{j}_{j+2}\},\cr
Q^{j+1}_{j+3} = \max\{lQ^{j}_{j+1},\,\, Q^{j}_{j+3}\},\cr
Q^{j+1}_{j+4} = \max\{lQ^{j}_{j+1},\,\, Q^{j}_{j+3},
Q^{j}_{j}+Q^{j}_{j+1}\}.
\end{cases}
$$

Rule $\Re_2$:

$$
\begin{cases}
Q^{j+1}_{j+1} = Q^{j}_{j+1},\cr
Q^{j+1}_{j+2} = Q^{j}_{j+2},\cr
Q^{j+1}_{j+3} = \max\{(2-l)Q^{j}_{j+1}+Q^{j}_{j},\,\,
Q^{j}_{j+3}\},\cr
Q^{j+1}_{j+4}
=\max\{(2-l)Q^{j}_{j+1}+Q^{j}_{j},\,\, Q^{j}_{j+3},\,\,
Q^{j}_{j}+Q^{j}_{j+1}\}.
\end{cases}
$$

Rule $\Re_3$:

$$
\begin{cases}
Q^{j+1}_{j+1} = Q^{j}_{j+1},\cr
 Q^{j+1}_{j+2} = Q^{j}_{j+2},\cr
Q^{j+1}_{j+3} = Q^{j}_{j+3},\cr
 Q^{j+1}_{j+4} =
\max\{Q^{j}_{j+3},\,\, Q^{j}_{j}+Q^{j}_{j+2}\}.
\end{cases}
$$

For   a sequence $I=(i_0, \ldots ,i_{r-1})$ we take a sequence of
rules $(\Re_{i_0}, \ldots ,\Re_{i_{r-1}})$ and construct a sequence
$\{Q_{j}(I, m)\}, j = 0,\ldots ,r+3$ in the following way.

For  $j = 0$ put
$$Q_0(m)=Q^0_0(m)=1,\,\,\, Q^0_t(m)=\theta^t_m,\,\,\, t=1,2,3.$$
For $j \ge 0$
given
$$(Q^j_j(m),\,\, Q^j_{j+1}(m),\,\, Q^j_{j+2}(m),\,\, Q^j_{j+3}(m))$$
we construct $$(Q^{j+1}_{j+1}(m),\,\, Q^{j+1}_{j+2}(m),\,\,
Q^{j+1}_{j+3}(m),\,\, Q^{j+1}_{j+4}(m))$$
 by the rule $\Re_{i_j}$
with $l=l_m.$

Now we define $Q_{j}(I,m)=Q^j_j(m)$ for $j\leq r$ and
$Q_{r+t,I}(m)=Q^r_{r+t}(m),\,\,\, t=1,2,3.$\\

The following proposition presents a result of computer calculation.

{\bf Proposition 2.}\,\,\, \emph{ Let $m=0, \ldots ,6$. For
any considered sequence of rules $I$ and defined sequence
$\{Q_{k}(I,m)\}$ one has}
$$(Q_{r+j}(I,m))^{\frac{1}{r+j}}\geq \theta_{m+1},\,\,\, j=0,1,2,3.$$
\\

Remind that the increasing sequence of remainders of best
approximations $\{q_k\}$ satisfies (\ref{1.2}) and Proposition 1.
When $l\in(0,2)$ all coefficients in inequalities (6), (7), (8) are
positive.

So we immediately deduce from Proposition 1 and Proposition 2 the
following statement:

{\bf Proposition 3.}\,\,\, \emph{ Suppose that
$$q_{i+j}\geq \lambda\theta_m^j,\,\,\, j=0,1,2,3,\,\,\,\, \lambda>0.$$
Then}
$$q_{r+j}\geq \lambda\theta^{r+j}_{m+1},\,\,\, j=0,1,2,3.$$\\

From (\ref{1.2}) it follows that for some positive $\lambda$ one has
$q_i\geq\lambda\theta_0^j.$ By Proposition 3 we see that
$q_{j+4r_0+2r_1}\geq\ \lambda\theta^{j+4r_0+2r_1}_{7}$ for any $j$.
 Theorem 3 is proved.
\\
\\

\vskip+1.0cm

Author's address:

Dept. of Number Theory

Fac. Mathematics and Mechanics

Moscow State University

119992 Moscow

Russia

e-mail: zzremi@gmail.com

\end{document}